\documentclass[10pt]{article}
\usepackage{cite}
\usepackage{mathrsfs}
\usepackage{amsfonts}
\usepackage{amsmath}
\usepackage{amsfonts,amssymb}
\usepackage{dsfont}
\usepackage{curves}
\usepackage{mathrsfs}
\usepackage{pifont}
\usepackage{amssymb}
\allowdisplaybreaks

\numberwithin{equation}{section}

\date{}

\def\BigRoman{\uppercase\expandafter{\romannumeral\number\count 255 }}
\def\Romannumeral{\afterassignment\BigRoman\count255=}

\begin{document}
\title{Signless Laplacian spectral radius for a $k$-extendable graph
}
\author{\small  Sizhong Zhou$^{1}$\footnote{Corresponding
author. E-mail address: zsz\_cumt@163.com (S. Zhou)}, Yuli Zhang$^{2}$\footnote{Corresponding
author. E-mail address: zhangyuli\_djtu@126.com (Y. Zhang)}\\
\small $1$. School of Science, Jiangsu University of Science and Technology,\\
\small Zhenjiang, Jiangsu 212100, China\\
\small $2$. School of Science, Dalian Jiaotong University,\\
\small Dalian, Liaoning 116028, China\\
}

\maketitle
\begin{abstract}
\noindent Let $k$ and $n$ be two nonnegative integers with $n\equiv0$ (mod 2), and let $G$ be a graph of order $n$ with a 1-factor. Then $G$ is 
said to be $k$-extendable for $0\leq k\leq\frac{n-2}{2}$ if every matching in $G$ of size $k$ can be extended to a 1-factor. In this paper, we
first establish a lower bound on the signless Laplacian spectral radius of $G$ to ensure that $G$ is $k$-extendable. Then we create some extremal
graphs to claim that all the bounds derived in this article are sharp.
\\
\begin{flushleft}

(2020) Mathematics Subject Classification: 05C50, 05C70
\end{flushleft}
\end{abstract}

\section{Introduction}

Graphs discussed in this paper are simple, undirected and connected. Let $G$ be a graph with vertex set $V(G)=\{v_1,v_2,\ldots,v_n\}$ and
edge set $E(G)$, where $|V(G)|=n$ and $\overline{G}$ be the complement of $G$. Denoted by $N_G(v)$ the neighbor set of the vertex $v$ in
$G$. The degree of the vertex $v$ is $d_G(v)=|N_G(v)|$. For $S\subseteq V(G)$, $G[S]$ denotes the subgraph of $G$ induced by $S$ and $G-S$
is the subgraph of $G$ induced by $V(G)\setminus S$. Given two vertex-disjoint graphs $G_1$ and $G_2$, the union of $G_1$ and $G_2$ is 
denoted by $G_1\cup G_2$ and the join of $G_1$ and $G_2$ is denoted by $G_1\vee G_2$. Let $K_n$ denote the complete graph of order $n$.

Let $A(G)$ denote the $(0,1)$-adjacency matrix of $G$ and $D(G)=diag(d_1,d_2,\ldots,d_n)$ denote the diagonal degree matrix of $G$, where
$d_i=d_G(v_i)$ for $1\leq i\leq n$. The signless Laplacian matrix $Q(G)$ of $G$ is defined as $Q(G)=D(G)+A(G)$. Obviously, $A(G)$ and $Q(G)$
are real symmetric matrices. The largest eigenvalues of $A(G)$ and $Q(G)$, denoted by $\rho(G)$ and $q(G)$, are called the spectral radius
and the signless Laplacian spectral radius of $G$, respectively.

For two positive integers $a$ and $b$ with $a\leq b$, a spanning subgraph $F$ of $G$ is called an $[a,b]$-factor if $a\leq d_G(v)\leq b$
for any $v\in V(G)$. If $a=b=1$, then an $[a,b]$-factor is a 1-factor (or a perfect matching). Let $G$ be a graph of order $n$ with a
1-factor. Then $G$ is said to be $k$-extendable for $0\leq k\leq\frac{n-2}{2}$ if every matching in $G$ of size $k$ can be extended to a
1-factor. In particular, $G$ is 0-extendable if and only if $G$ contains a 1-factor.

Many researchers have attempted to find sufficient conditions for the existence of 1-factors by utilizing various graphic parameters.
Tutte \cite{T} obtained a characterization for a graph with a 1-factor. Anderson \cite{A1,A2} investigated the relationships between
binding numbers and 1-factors in graphs and presented two binding number conditions for the existence of 1-factors in graphs. Sumner
\cite{S} showed a sufficient condition for a graph to possess a 1-factor. Niessen \cite{N} provided a neighborhood union condition for
the existence of 1-factors in graphs. Enomoto \cite{E} derived a toughness condition for a graph to admit a 1-factor. Plummer \cite{P1}
first introduced the concept of $k$-extendable graph and posed some properties of $k$-extendable graphs. Up to now, much attention has
been paid on various graphic parameters of $k$-extendable graphs, such as binding number \cite{C,RW}, connectivity \cite{LY,P2}, minimum
degree \cite{AC}, independence number \cite{MV,AC2} and genus \cite{P3}. Much effort has been devoted to finding sufficient conditions
for the existence of $[1,2]$-factors (see \cite{Ke,KKK,ZSL,ZWB,ZWX,Zhr,Zd,ZB,GW,Ls,WZi,Wp}) and $[a,b]$-factors (see
\cite{KL,M,WZo,WZr,GWC,ZL,Za1,Za2,Za3,Zr,ZLX,ZXS}) in graphs.

The main goal of this paper is to study the existence of $k$-extendable graphs from a spectral perspective. Recall that $G$ is 0-extendable
if and only if $G$ has a 1-factor. In the past few years, lots of researchers focused on finding the connections between the spectral
radius and 1-factors in graphs. O \cite{Os} provided a spectral radius condition to guarantee that a connected graph has a 1-factor. By
imposing the minimum degree of a graph as a parameter, Liu, Liu and Feng \cite{LLF} extended O's result \cite{Os} in a connected graph.
In this paper, we study the existence of $k$-extendable graphs and obtain a signless Laplacian spectral radius condition for graphs to be 
$k$-extendable.

\medskip

\noindent{\textbf{Theorem 1.1.}} Let $k$ and $n$ be two nonnegative integers with $n\equiv0$ (mod 2), and let $G\neq K_{2k}\vee(K_{n-2k-1}\cup K_1)$ 
be a connected graph of even order $n$ with $n\geq2k+4$.\\
(\romannumeral1) For $n\notin\{2k+6,2k+8\}$, if $q(G)>\theta(k,n)$, then $G$ is $k$-extendable, where $\theta(k,n)$ is the largest root of
$x^{3}-(3n+2k-7)x^{2}+(2n^{2}+6kn-7n-24k)x-2(2k+1)(n-3)(n-4)=0$.\\
(\romannumeral2) For $n=2k+6$, if $q(G)>3k+4+\sqrt{k^{2}+12k+12}$, then $G$ is $k$-extendable.\\
(\romannumeral3) For $n=2k+8$, if $q(G)>3k+6+\sqrt{k^{2}+16k+24}$, then $G$ is $k$-extendable.

\medskip

The proof of Theorem 1.1 will be provided in Section 3.

\section{Preliminary lemmas}

In this section, we put forward some necessary preliminary lemmas, which are very important to the proofs of our main results.

Chen \cite{C} established a necessary and sufficient condition for the existence of $k$-extendable graphs.

\medskip

\noindent{\textbf{Lemma 2.1}} (\cite{C}). Let $k\geq1$ be an integer. Then a graph $G$ is $k$-extendable if and only if
$$
o(G-S)\leq|S|-2k
$$
for any $S\subseteq V(G)$ such that $G[S]$ contains $k$ independent edges, where $o(G-S)$ denotes the number of odd components in $G-S$.

\medskip

\noindent{\textbf{Lemma 2.2}} (\cite{SYZL}). Let $G$ be a connected graph. If $H$ is a subgraph of $G$, then $q(H)\leq q(G)$.
If $H$ is a proper subgraph of $G$, then $q(H)<q(G)$.

\medskip

\noindent{\textbf{Lemma 2.3}} (\cite{CS}). Let $n\geq2$ be an integer, and $K_n$ be a complete graph of order $n$. Then $q(K_n)=2n-2$.

\medskip

In what follows, we explain the concepts of equitable matrices and equitable partitions.

\medskip

\noindent{\textbf{Definition 2.4}} (\cite{BH1}). Let $M$ be a real matrix of order $n$ described in the following block form
\begin{align*}
M=\left(
  \begin{array}{ccc}
    M_{11} & \cdots & M_{1r}\\
    \vdots & \ddots & \vdots\\
    M_{r1} & \cdots & M_{rr}\\
  \end{array}
\right),
\end{align*}
where the blocks $M_{ij}$ are $n_i\times n_j$ matrices for any $1\leq i,j\leq r$ and $n=n_1+n_2+\cdots+n_r$. For $1\leq i,j\leq r$, let $b_{ij}$
denote the average row sum of $M_{ij}$, that is, $b_{ij}$ is the sum of all entries in $M_{ij}$ divided by the number of rows. Then $B(M)=(b_{ij})$
(simply by $B$) is called a quotient matrix of $M$. If for every pair $i,j$, $M_{ij}$ admits constant row sum, then $B$ is called an equitable
quotient matrix of $M$ and the partition is called equitable.

\medskip

\noindent{\textbf{Lemma 2.5}} (\cite{YYSX}). Let $B$ be an equitable matrix of $M$ as defined in Definition 2.4, and $M$ be
a nonnegative matrix. Then $\rho_1(B)=\rho_1(M)$, where $\rho_1(B)$ and $\rho_1(M)$ denote the largest eigenvalues of the matrices $B$ and $M$.

\section{The proof of Theorem 1.1}

In this section, we prove Theorem 1.1, which provides a sufficient condition via the signless Laplacian spectral radius of a connected graph
to ensure that the graph is $k$-extendable.

\medskip

\medskip

\noindent{\it Proof of Theorem 1.1.} Suppose, to the contrary, that $G$ is not $k$-extendable, according to Lemma 2.1, there exists some nonempty
subset $S$ of $V(G)$ such that $|S|\geq2k$ and $o(G-S)>|S|-2k$. Since $n$ is even, $o(G-S)$ and $|S|$ possess the same parity. Thus, we deduce
$$
o(G-S)\geq|S|-2k+2.
$$
Select such a connected graph $G$ of order $n$ so that its signless Laplacian spectral radius is as large as possible.

Together with Lemma 2.2 and the choice of $G$, the induced subgraph $G[S]$ and every connected component of $G-S$ are complete graphs, respectively.
Furthermore, all components of $G-S$ are odd and $G$ is just the graph $G[S]\vee(G-S)$.

For convenience, let $o(G-S)=q$ and $|S|=s$. Then $q\geq s-2k+2$. Assume that $G_1,G_2,\ldots,G_q$ are all the components of $G-S$ with
$n_i=|V(G_i)|$ and $n_1\geq n_2\geq\cdots\geq n_q$. Then $G=K_s\vee(K_{n_1}\cup K_{n_2}\cup\cdots\cup K_{n_q})$.

\medskip

\noindent{\bf Claim 1.} $n_2=n_3=\cdots=n_q=1$.

\noindent{\it Proof.} If $n_2\geq3$, then we let $G'=K_s\vee(K_{n_1+2}\cup K_{n_2-2}\cup K_{n_3}\cup\cdots\cup K_{n_q})$. Note that
$o(G'-S)=o(G-S)=q\geq s-2k+2$. Denote the vertex set of $G$ by $V(G)=V(K_s)\cup V(K_{n_1})\cup V(K_{n_2})\cup\cdots\cup V(K_{n_q})$. Let $Y$
be the Perron vector of $Q(G)$, and let $Y(v)$ be the entry of $Y$ corresponding to the vertex $v\in V(G)$. By symmetry, it is obvious that
all vertices of $K_s$ (resp. $K_{n_1},K_{n_2},\cdots,K_{n_q}$) have the same entries in $Y$. Hence, we can suppose $Y(v_0)=y_0$ for every
$v_0\in V(K_s)$, $Y(v_1)=y_1$ for every $v_1\in V(K_{n_1})$, $Y(v_2)=y_2$ for every $v_2\in V(K_{n_2}), \cdots, Y(v_q)=y_q$ for every
$v_q\in V(K_{n_q})$. Then
\begin{equation}\label{eq:3.1}
\begin{cases}
q(G)y_1=sy_0+(s+2n_1-2)y_1,\\
q(G)y_2=sy_0+(s+2n_2-2)y_2.
\end{cases}
\end{equation}
It follows from \eqref{eq:3.1} that
\begin{align}\label{eq:3.2}
(q(G)-s-2n_1+2)y_1=(q(G)-s-2n_2+2)y_2.
\end{align}
Note that $K_{s+n_1}$ and $K_{s+n_2}$ are two proper subgraphs of $G$. Using Lemmas 2.2 and 2.3, we get
\begin{align*}
q(G)>&\max\{q(K_{s+n_1}),q(K_{s+n_2})\}\\
=&\max\{2(s+n_1)-2,2(s+n_2)-2\}\\
>&\max\{s+2n_1-2,s+2n_2-2\}.
\end{align*}
Together with \eqref{eq:3.2} and $n_1\geq n_2$, we infer $y_1\geq y_2$. According to the Rayleigh quotient, we derive
\begin{align*}
q(G')-q(G)\geq&Y^{T}(Q(G')-Q(G))Y\\
=&2n_1y_1(y_1+y_2)+2n_1y_2(y_1+y_2)-8(n_2-2)y_2^{2}\\
\geq&8n_1y_2^{2}-8(n_2-2)y_2^{2}\\
=&8y_2^{2}(n_1-n_2+2)\\
>&0.
\end{align*}
Hence, $q(G')>q(G)$, which is a contradiction to the choice of $G$. Thus, we deduce $n_2=1$.

Recall that $n_2\geq n_3\geq\cdots\geq n_q\geq1$. Combining this with $n_2=1$, we infer $n_2=n_3=\cdots=n_q=1$. Claim 1 is proved. \hfill $\Box$

In what follows, we are to verify $q=s-2k+2$. Note that $q\geq s-2k+2$, and $q$ and $s$ have the same parity. Consequently, we can suppose
$q\geq s-2k+4$. We construct a new graph $G''=K_s\vee(K_{n_1+2}\cup(q-3)K_1)$. Clearly, $G$ is a proper subgraph of $G''$ and
$o(G''-S)=o(G-S)-2=q-2\geq s-2k+2$. Together with Lemma 2.2, $q(G'')>q(G)$, which is a contradiction to the choice of $G$. Thus, we infer
$q\leq s-2k+2$. On the other hand, $q\geq s-2k+2$. Hence, we obtain
\begin{align}\label{eq:3.3}
q=s-2k+2.
\end{align}

By virtue of \eqref{eq:3.3}, Claim 1, $n=s+n_1+n_2+\cdots+n_q$ and $G=K_s\vee(K_{n_1}\cup K_{n_2}\cup\cdots\cup K_{n_q})$, we have
$G=K_s\vee(K_{n_1}\cup(q-1)K_1)=K_s\vee(K_{n_1}\cup(s-2k+1)K_1)$ and $n_1=n-s-(q-1)=n-2s+2k-1$. If $s=2k$, then $G=K_{2k}\vee(K_{n-2k-1}\cup K_1)$,
which is a contradiction to the condition of this theorem. Hence, $s\geq2k+1$. The following proof will be divided into two cases by the value of
$n_1$.

\noindent{\bf Case 1.} $n_1\geq3$.

In this case, $n=n_1+2s-2k+1\geq2s-2k+4$. Recall that $G=K_s\vee(K_{n_1}\cup(s-2k+1)K_1)$ and $n_1=n-2s+2k-1$. Consider the partition
$V(G)=V(K_s)\cup V(K_{n_1})\cup V((s-2k+1)K_1)$. The corresponding quotient matrix of $Q(G)$ equals
\begin{align*}
B_1=\left(
  \begin{array}{ccc}
    n+s-2 & n-2s+2k-1 & s-2k+1\\
    s & 2n-3s+4k-4 & 0\\
    s & 0 & s\\
  \end{array}
\right).
\end{align*}
Then the characteristic polynomial of $B_1$ is
\begin{align*}
f_1(x)=&x^{3}-(3n-s+4k-6)x^{2}+(2n^{2}+sn+4kn-8n-4s^{2}-4s+8ks-8k+8)x\\
&-2sn^{2}+4s^{2}n-8ksn+10sn-2s^{3}+8ks^{2}-10s^{2}-8k^{2}s+20ks-12s.
\end{align*}
In view of Lemma 2.5, the largest root, say $q_1$, of $f_1(x)=0$ equals the signless Laplacian spectral radius of $G$. Consequently, we possess
$f_1(q_1)=0$ and $q(G)=q_1$.

Note that $K_s\vee(n-s)K_1$ is a proper subgraph of $G$. From Lemma 2.2, we infer $q_1=q(G)>q(K_s\vee(n-s)K_1)$. Consider the partition
$V(K_s\vee(n-s)K_1)=V(K_s)\cup V((n-s)K_1)$. The corresponding quotient matrix of $Q(K_s\vee(n-s)K_1)$ has the following form
\begin{align*}
B_2=\left(
  \begin{array}{ccc}
    n+s-2 & n-s\\
    s & s\\
  \end{array}
\right).
\end{align*}
Then the characteristic polynomial of $B_2$ equals
\begin{align*}
f_2(x)=x^{2}-(n+2s-2)x+2s(s-1).
\end{align*}
In terms of Lemma 2.5, the largest root, say $q_2$, of $f_2(x)=0$ equals $q(K_s\vee(n-s)K_1)$. And so
\begin{align}\label{eq:3.4}
q(K_s\vee(n-s)K_1)=q_2=\frac{n+2s-2+\sqrt{(n+2s-2)^{2}-8s(s-1)}}{2}.
\end{align}
Together with $q_1=q(G)>q(K_s\vee(n-s)K_1)$, we get
\begin{align}\label{eq:3.5}
q_1>q_2=\frac{n+2s-2+\sqrt{(n+2s-2)^{2}-8s(s-1)}}{2}.
\end{align}
Note that $f_1(q_1)=0$. By a direct calculation, we have
\begin{align}\label{eq:3.6}
\varphi(q_1)=\varphi(q_1)-f_1(q_1)=(s-2k-1)g_1(q_1),
\end{align}
where $g_1(q_1)=-q_1^{2}+(-n+4s+8)q_1+2n^{2}-4sn-14n+2s^{2}-4ks+12s+24$. Utilizing \eqref{eq:3.5} and $n\geq2s-2k+4\geq s+5$, we derive
$$
-\frac{-n+4s+8}{2\times(-1)}<n+s-2<\frac{n+2s-2+\sqrt{(n+2s-2)^{2}-8s(s-1)}}{2}<q_1.
$$
Consequently, we deduce
\begin{align}\label{eq:3.7}
g_1(q_1)<&g_1\left(\frac{n+2s-2+\sqrt{(n+2s-2)^{2}-8s(s-1)}}{2}\right)\nonumber\\
=&n^{2}-5sn-7n+6s^{2}-4ks+18s+14\nonumber\\
&-(n-s-5)\sqrt{(n+2s-2)^{2}-8s(s-1)}\nonumber\\
\leq&n^{2}-5sn-7n+6s^{2}-4ks+18s+14-n(n-s-5)\nonumber\\
=&-4sn-2n+6s^{2}-4ks+18s+14\nonumber\\
\leq&-4s(2s-2k+4)-2(2s-2k+4)+6s^{2}-4ks+18s+14\nonumber\\
=&-2s^{2}+4ks-2s+4k+6.
\end{align}

For $s\geq2k+2$, it follows from \eqref{eq:3.7} that
\begin{align}\label{eq:3.8}
g_1(q_1)<&-2s^{2}+4ks-2s+4k+6\nonumber\\
\leq&-2s(2k+2)+4ks-2s+4k+6\nonumber\\
=&-6s+4k+6\nonumber\\
<&0.
\end{align}

According to \eqref{eq:3.6}, \eqref{eq:3.8} and $s\geq2k+1$, we infer
$$
\varphi(q_1)=(s-2k-1)g_1(q_1)\leq0,
$$
which yields $q(G)=q_1\leq\theta(k,n)$, which is a contradiction to $q(G)>\theta(k,n)$.

\noindent{\bf Case 2.} $n_1=1$.

In this case, we possess $G=K_s\vee(s-2k+2)K_1=K_s\vee(n-s)K_1$ and $n=2s-2k+2$. By virtue of \eqref{eq:3.4}, we obtain
$$
q(G)=q(K_s\vee(n-s)K_1)=q_2=\frac{n+2s-2+\sqrt{(n+2s-2)^{2}-8s(s-1)}}{2}.
$$
Note that $f_2(q_2)=0$. By a direct computation, we possess
\begin{align}\label{eq:3.9}
\varphi(q_2)=&\varphi(q_2)-q_2f_2(q_2)\nonumber\\
=&-(2n-2s+2k-5)q_2^{2}+(2n^{2}+6kn-7n-24k-2s^{2}+2s)q_2\nonumber\\
&-2(2k+1)n^{2}+14(2k+1)n-24(2k+1)\nonumber\\
=&-(s+2k+1)n^{2}+(3s+2k+4)sn+7(2k+1)n-2s^{3}-(20k+8)s-14(2k+1)\nonumber\\
&+(-sn+2kn+n+s^{2}-2ks+4s-10k-5)\sqrt{(n+2s-2)^{2}-8s(s-1)}.
\end{align}

If $s=2k+1$, then $n=2k+4$ and $q(G)=3k+2+\sqrt{k^{2}+8k+4}=\theta(k,2k+4)$, which contradicts $q(G)>\theta(k,n)$ for $n=2k+4$. If $s=2k+2$, then
$n=2k+6$ and $q(G)=3k+4+\sqrt{k^{2}+12k+12}$, which contradicts $q(G)>3k+4+\sqrt{k^{2}+12k+12}$ for $n=2k+6$. If $s=2k+3$, then $n=2k+8$ and
$q(G)=3k+6+\sqrt{k^{2}+16k+24}$, which contradicts $q(G)>3k+6+\sqrt{k^{2}+16k+24}$ for $n=2k+8$. In what follows, we consider $s\geq2k+4$.

Recall that $n=2s-2k+2$. According to $s\geq2k+4$, we easily see
\begin{align*}
(n+2s-2)^{2}-8s(s-1)=&8s^{2}-8(2k-1)s+4k^{2}\nonumber\\
\geq&4s^{2}+4(2k+4)s-8(2k-1)s+4k^{2}\nonumber\\
=&4s^{2}-8(k-3)s+4k^{2}\nonumber\\
\geq&4s^{2}-8(k-2)s+8(2k+4)+4k^{2}\nonumber\\
=&(2s-2k+4)^{2}+32k+16\nonumber\\
>&(2s-2k+4)^{2}\nonumber\\
=&(n+2)^{2}
\end{align*}
and
$$
-sn+2kn+n+s^{2}-2ks+4s-10k-5=-(s-2k-1)(s-2k-3)<0.
$$
Combining these with \eqref{eq:3.9}, $s\geq2k+4$ and $n=2s-2k+2$, we deduce
\begin{align}\label{eq:3.10}
\varphi(q_2)=&-(s+2k+1)n^{2}+(3s+2k+4)sn+7(2k+1)n-2s^{3}-(20k+8)s-14(2k+1)\nonumber\\
&+(-sn+2kn+n+s^{2}-2ks+4s-10k-5)\sqrt{(n+2s-2)^{2}-8s(s-1)}\nonumber\\
<&-(s+2k+1)n^{2}+(3s+2k+4)sn+7(2k+1)n-2s^{3}-(20k+8)s-14(2k+1)\nonumber\\
&-(s-2k-1)(s-2k-3)(n+2)\nonumber\\
=&-2s^{3}+(8k+6)s^{2}-(8k^{2}+4k-12)s-(2k+1)(8k+16)\nonumber\\
:=&p(s).
\end{align}

We easily see
$$
p'(s)=-6s^{2}+2(8k+6)s-8k^{2}-4k+12
$$
and
$$
p''(s)=-12s+2(8k+6).
$$
Recall that $s\geq2k+4$. Then $p''(s)=-12s+2(8k+6)\leq-12(2k+4)+2(8k+6)=-8k-36<0$, which implies that $p'(s)$ is decreasing in the interval
$[2k+4,+\infty)$. Thus, $p'(s)\leq p'(2k+4)=-6(2k+4)^{2}+2(8k+6)(2k+4)-8k^{2}-4k+12=-12k-36<0$, which yields that $p(s)$ is decreasing in the
interval $[2k+4,+\infty)$. Thus, $p(s)\leq p(2k+4)=-2(2k+4)^{3}+(8k+6)(2k+4)^{2}-(8k^{2}+4k-12)(2k+4)-(2k+1)(8k+16)=0$. Together with \eqref{eq:3.10},
we infer $\varphi(q_2)<p(s)\leq0$, which implies $q(G)=q_2<\theta(k,n)$, a contradiction to the condition. This completes the proof of
Theorem 1.1. \hfill $\Box$

\section{Concluding remark}

In this section, we claim that the bounds derived in Theorem 1.1 are best possible.

\medskip

\noindent{\textbf{Theorem 4.1.}} Let $k$ and $n$ be two nonnegative integers with $n\equiv0$ (mod 2), and let $\theta(k,n)$ be the largest root
of $x^{3}-(3n+2k-7)x^{2}+(2n^{2}+6kn-7n-24k)x-2(2k+1)(n-3)(n-4)=0$. For $n\geq2k+4$ and $n\notin\{2k+6,2k+8\}$, we have 
$q(K_{2k+1}\vee(K_{n-2k-3}\cup2K_1))=\theta(k,n)$ and $K_{2k+1}\vee(K_{n-2k-3}\cup2K_1)$ is not $k$-extendable. For $n=2k+6$, we possess 
$q(K_{2k+2}\vee4K_1)=3k+4+\sqrt{k^{2}+12k+12}$ and $K_{2k+2}\vee4K_1$ is not $k$-extendable. For $n=2k+8$, we admit 
$q(K_{2k+3}\vee5K_1)=3k+6+\sqrt{k^{2}+16k+24}$ and $K_{2k+3}\vee5K_1$ is not $k$-extendable.

\medskip

\noindent{\it Proof.} Consider the partition $V(K_{2k+1}\vee(K_{n-2k-3}\cup2K_1))=V(K_{2k+1})\cup V(K_{n-2k-3})\cup V(2K_1)$. The corresponding 
quotient matrix of $Q(K_{2k+1}\vee(K_{n-2k-3}\cup2K_1))$ is equal to
\begin{align*}
B_1=\left(
  \begin{array}{ccc}
    n+2k-1 & n-2k-3 & 2\\
    2k+1 & 2n-2k-7 & 0\\
    2k+1 & 0 & 2k+1\\
  \end{array}
\right).
\end{align*}
Then the characteristic polynomial of the matrix $B_1$ is equal to $x^{3}-(3n+2k-7)x^{2}+(2n^{2}+6kn-7n-24k)x-2(2k+1)(n-3)(n-4)$. In terms of
Lemma 2.5, the largest root $\theta(k,n)$ of $x^{3}-(3n+2k-7)x^{2}+(2n^{2}+6kn-7n-24k)x-2(2k+1)(n-3)(n-4)=0$ equals $q(K_{2k+1}\vee(K_{n-2k-3}\cup2K_1))$. 
Namely, $q(K_{2k+1}\vee(K_{n-2k-3}\cup2K_1))=\theta(k,n)$. Write $S=V(K_{2k+1})$. Then $o(K_{2k+1}\vee(K_{n-2k-3}\cup2K_1)-S)=3>1=(2k+1)-2k=|S|-2k$. 
By virtue of Lemma 2.1, the graph $K_{2k+1}\vee(K_{n-2k-3}\cup2K_1)$ is not $k$-extendable.

\medskip

Consider the partition $V(K_{2k+2}\vee4K_1)=V(K_{2k+2})\cup V(4K_1)$. The corresponding quotient matrix of $Q(K_{2k+2}\vee4K_1)$ equals
\begin{align*}
B_2=\left(
  \begin{array}{ccc}
    4k+6 & 4\\
    2k+2 & 2k+2\\
  \end{array}
\right).
\end{align*}
Then the characteristic polynomial of the matrix $B_2$ is $x^{2}-(6k+8)x+(2k+2)(4k+2)$. It follows from Lemma 2.5 that the largest root of 
$x^{2}-(6k+8)x+(2k+2)(4k+2)=0$ equals $q(K_{2k+2}\vee4K_1)$. Thus, we possess $q(K_{2k+2}\vee4K_1)=3k+4+\sqrt{k^{2}+12k+12}$. Let $S=V(K_{2k+2})$. 
Then $o(K_{2k+2}\vee4K_1-S)=4>2=(2k+2)-2k=|S|-2k$. Applying Lemma 2.1, the graph $K_{2k+2}\vee4K_1$ is not $k$-extendable.

\medskip

Consider the partition $V(K_{2k+3}\vee5K_1)=V(K_{2k+3})\cup V(5K_1)$. The corresponding quotient matrix of $Q(K_{2k+3}\vee5K_1)$ is equal to
\begin{align*}
B_3=\left(
  \begin{array}{ccc}
    4k+9 & 5\\
    2k+3 & 2k+3\\
  \end{array}
\right).
\end{align*}
Then the characteristic polynomial of the matrix $B_3$ equals $x^{2}-(6k+12)x+(2k+3)(4k+4)$. Utilizing Lemma 2.5, the largest root of 
$x^{2}-(6k+12)x+(2k+3)(4k+4)=0$ is equal to $q(K_{2k+3}\vee5K_1)$. Thus, we deduce $q(K_{2k+3}\vee5K_1)=3k+6+\sqrt{k^{2}+16k+24}$. Write 
$S=V(K_{2k+3})$. Then $o(K_{2k+3}\vee5K_1-S)=5>3=(2k+3)-2k=|S|-2k$. It follows from Lemma 2.1 that the graph $K_{2k+3}\vee5K_1$ is not 
$k$-extendable. \hfill $\Box$

\medskip

\section*{Data availability statement}

My manuscript has no associated data.

\section*{Declaration of competing interest}

The authors declare that they have no conflicts of interest to this work.



\end{document}